\documentclass[a4paper, 12pt, roman]{article}
\usepackage{MRHStyle}
\usepackage{JournalOfExpMath}

\usepackage[square]{natbib}
\setcitestyle{citesep={;}}

\renewcommand{\phi}{\ensuremath{\varphi}}

\begin{document}
\title{\Large \textbf{An extension of the Geometric Modulus Principle to holomorphic and harmonic functions}}
\author{\Large Matt Hohertz \\ Department of Mathematics, Rutgers University \\ matt.hohertz@rutgers.edu}

\maketitle

\begin{abstract}
	Kalantari's Geometric Modulus Principle describes the local behavior of the modulus of a polynomial.  Specifically, if $p(z) = a_0 + \sum_{j=k}^n a_j\inps{z-z_0}^j,\;a_0a_ka_n \neq 0$, then the complex plane near $z = z_0$ comprises $2k$ sectors of angle $\frac{\pi}{k}$, alternating between \ic{arguments of ascent} (angles $\theta$ where $\aval{p(z_0 + te^{i\theta})} > \aval{p(z_0)}$ for small $t$) and \ic{arguments of descent} (where the opposite strict inequality holds).  In this paper, we generalize the Geometric Modulus Principle to holomorphic and harmonic functions.  As in Kalantari's original paper, we use these extensions to give succinct, elegant new proofs of some classical theorems from analysis.
\end{abstract}

\section*{Introduction}
	For a holomorphic function $p$ and fixed $z_0\in\C$, we define a value $\theta\in[0,2\pi)$ to be 
	\begin{enumerate}
		\item an \ic{argument of ascent} if there exists $t^* > 0$ such that 
			\begin{equation}
				\aval{p(z_0 + te^{i\theta})} > \aval{p(z_0)},\quad t\in\inps{0, t^*}; \mbox{ and}
			\end{equation}
		\item an \ic{argument of descent} if there exists $t^* > 0$ such that 
			\begin{equation}
				\aval{p(z_0 + te^{i\theta})} < \aval{p(z_0)},\quad t\in\inps{0, t^*}.
			\end{equation}
	\end{enumerate}  Moreover, we define the \ic{cone of ascent} $C_{a}(p, z_0)$ to be the intersection of the set of arguments of ascent with $[0, 2\pi)$ and define the \ic{cone of descent} $C_d(p, z_0)$ analogously. 
	\par
	Kalantari's Geometric Modulus Principle is the following.
	\begin{theorem}[Geometric Modulus Principle\footnotemark] \label{thm:gmp}
		Let $p(z)$ be a nonconstant polynomial. If $p(z_0) = 0$ then $C_a(p, z_0) = [0, 2\pi)$. If $p(z_0)\neq 0$ then $C_a(p, z_0)$ and $C_d(p, z_0)$ partition the unit circle into alternating sectors of angle $\frac{\pi}{k}$. Specifically, if $re^{i\alpha}$ is the polar form of $\overline{p(z_0)}p^{(k)}(z_0)$ then $\theta$ is an argument of ascent if
		\begin{equation} \label{eq:aor}
			\frac{-2\alpha - \pi}{2k} < \theta < \frac{-2\alpha + \pi}{2k}\quad \inps{\mbox{mod }\frac{2\pi}{k}}.
		\end{equation}
		\footnotetext{\cite[Theorem 1]{geomod}}
The interiors of the remaining sectors are arguments of descent.  A boundary point between two sectors is either an argument of ascent or an argument of descent.
\end{theorem}
Except in its last sentence, Theorem \ref{thm:gmp} is a straightforward corollary of the following lemma.
\begin{lemma} \label{lemma:cos}
	Let $p(z)$ be a polynomial of degree $n\geq 1$. Suppose that $z_0$ is a complex number, $p(z_0)\neq 0$, and $k:=\min_{\ell\geq 1}\:\sbld{\ell : p^{(\ell)}(z_0)\neq 0}$. Given a real
number $\theta$, define $F_{\theta}(t):=\aval{p\inps{z_0 + te^{i\theta}}}^2$.  Finally, let $re^{i\alpha}$ be the polar form of $\overline{p(z_0)}p^{(k)}\inps{z_0}$.  Then $$F_{\theta}^{(k)}\inps{0} = 2r\cos\inps{\alpha + k\theta}.$$
	\begin{proof}
		We paraphrase the proof from \cite{geomod}.  First, by Taylor's Theorem we assume without loss of generality that $z_0 = 0$.  Then
		\begin{align}
			F_{\theta}(t) &= p(te^{i\theta})\cdot\overline{p(te^{i\theta})}\\
			&= \inps{a_0+\sum_{j\geq k} a_j t^j e^{ij\theta}}\cdot\inps{\overline{a_0}+\sum_{j\geq k} \overline{a_j} t^j e^{-ij\theta}}  \label{eq:realanal} \\
			&= \aval{a_0}^2 + t^k\cdot \inps{\overline{a_0}a_ke^{ik\theta} + a_0\overline{a_k}e^{-ik\theta}}  + t^{k+1}\cdot \inps{\cdots}. \label{eq:cprod}
		\intertext{The terms of this Maclaurin series correspond to the derivatives of $F_{\theta}(t)$ at $t = 0$, so that}
			F_{\theta}^{(k)}\inps{0} &= k!\cdot\inps{\overline{a_0}a_ke^{ik\theta} + a_0\overline{a_k}e^{-ik\theta}}.  \label{eq:match} 
		\intertext{Continuing, we obtain}
			F_{\theta}^{(k)}\inps{0} &= 2\cdot Re\inps{\overline{p(0)}p^{(k)}(0)e^{ik\theta}} \\
			&= 2\cdot Re\inps{re^{i\inps{\alpha + k\theta}}} \\
			&= 2r\cos\inps{\alpha + k\theta}.
		\end{align}
	\end{proof}
\end{lemma}
Thus the proof of Theorem \ref{thm:gmp} is very simple.
\begin{proof}[Theorem \ref{thm:gmp}]
	At $\theta$ such that $F_{\theta}^{(k)}\inps{0}\neq 0$, the behavior of $\aval{p(z_0 + te^{i\theta})}$ for small positive $t$ follows directly from Lemma \ref{lemma:cos}.  On the other hand, $F_{\theta}(t)$ is real analytic and not identically zero, so if $F_{\theta}^{(k)}\inps{0} = 0$ then there exists $m > k$ such that $F_{\theta}^{(m)}\inps{0}\neq 0$.  The sign of this derivative determines whether $\theta$ is an argument of ascent or descent.
\end{proof}
\section{Extension to holomorphic functions}
Extending the Geometric Modulus Principle requires the following lemma.
\begin{lemma} \label{lemma:egmp}
	In Lemma \ref{lemma:cos}, the requirement that $p(z)$ be a non-constant polynomial may be relaxed to a requirement that $p(z)$ be non-constant and holomorphic near $z_0$.
	\begin{proof}
		The proof of Lemma \ref{lemma:cos} requires only that $p(z)$ be complex analytic (equivalently, holomorphic) with at least one non-constant term.  In particular:
		\begin{enumerate}[label=(\roman*)]
			\item By Mertens' Theorem, Equation \eqref{eq:cprod} continues to hold as long as one of the two series converges absolutely.  But holomorphic functions converge absolutely within their discs of convergence.
			\item Equation \eqref{eq:match} continues to hold because the series of Equation \eqref{eq:cprod} can be differentiated term-by-term in its interval of convergence \cite[Theorem~3]{apostol_1952}.
		\end{enumerate}
	\end{proof}
\end{lemma}

\begin{figure}[h]
    		\centering
    		\includegraphics[width=0.75\textwidth]{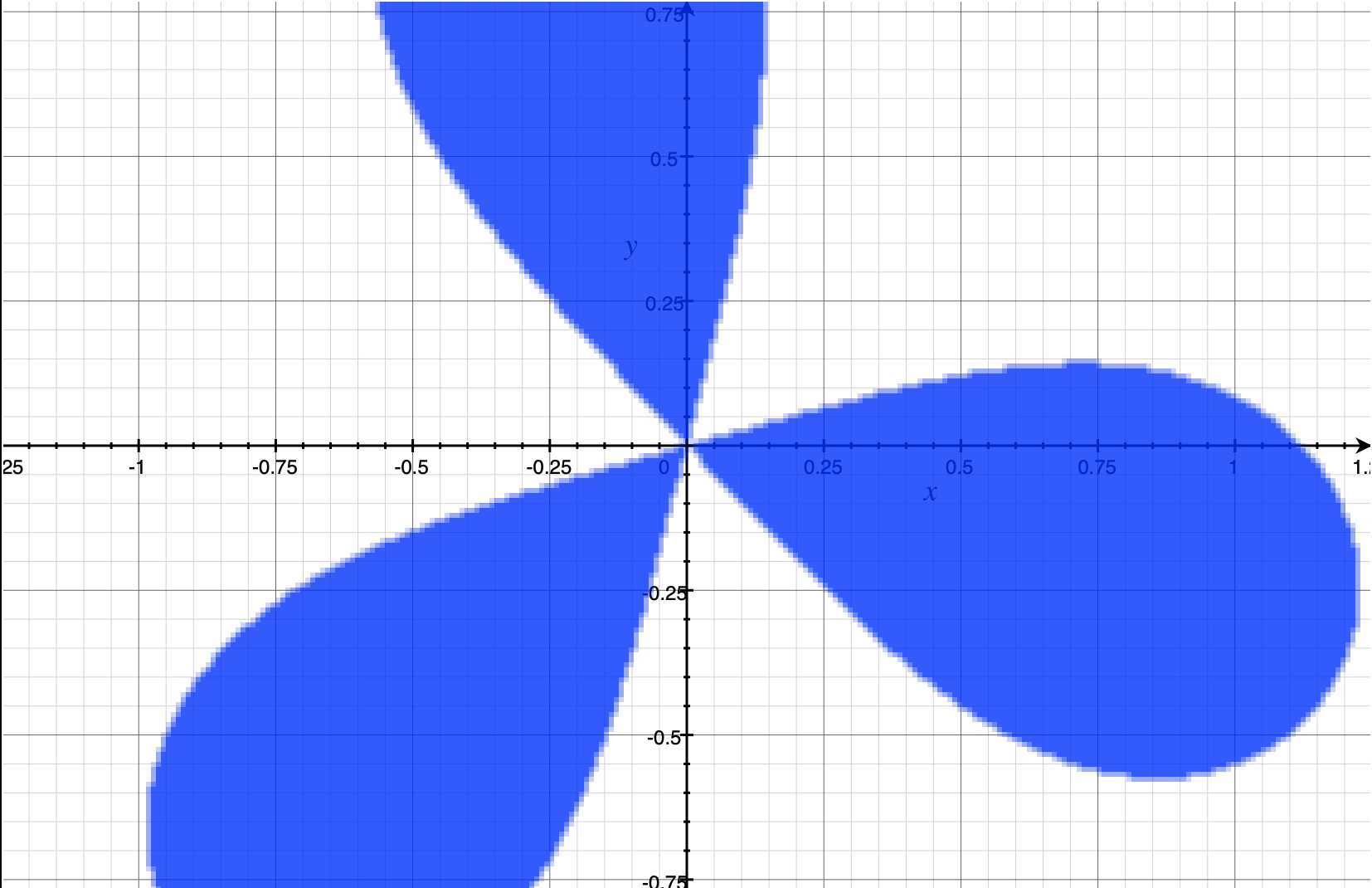}
  		  \caption{The lemniscate $\aval{\frac{1}{1-e^{i\pi/4}z^3}} > 1$.}
  		  \label{fig:GeoMod1}
\end{figure}

\begin{theorem}[Extended Geometric Modulus Principle] \label{thm:egmp}
	In Theorem \ref{thm:gmp}, the requirement that $p(z)$ be a non-constant polynomial may be relaxed to a requirement that $p(z)$ be non-constant and holomorphic near $z_0$.
	\begin{proof}
		The proof is the same as that of Theorem \ref{thm:gmp}, except with Lemma \ref{lemma:egmp} in place of Lemma \ref{lemma:cos}.
	\end{proof}
\end{theorem}
\paragraph{Example} Define $\omega:=e^{i\frac{\pi}{4}}$ and $f(z):=\frac{1}{1-\omega z^3}$.  Near $z_0 = 0$, $f(z)$ has the expansion
\begin{equation}
	f(z) = 1 + \omega z^3 + iz^6 + \cdots,
\end{equation}
so that
\begin{align}
	re^{i\alpha} &= \overline{f(0)}\cdot f^{(3)}(0) \\
	&= 6 \omega,
\end{align} yielding $\alpha = \frac{\pi}{4}$.  Thus we expect that
\begin{equation}
	C_a(f, 0) = \inps{-\frac{\pi}{4}, \frac{\pi}{12}} \cup \inps{\frac{5\pi}{12}, \frac{3\pi}{4}} \cup \inps{\frac{13\pi}{12}, \frac{17\pi}{12}}\quad \inps{\mbox{mod }2\pi}
\end{equation} up to boundary points.  Figure \ref{fig:GeoMod1}, which graphs the lemniscate $\sbld{\aval{f(z)} > \aval{f(0)}}$, bears this prediction out.
 \par

\section{Extension to harmonic functions}
In this section we extend the Geometric Modulus Principle to harmonic functions, using the fact that all harmonic functions are locally the real part of some holomorphic function.  Define an \ic{argument of ascent} (resp., \ic{argument of descent}) of a harmonic function $u$ at $z_0$ to be an angle $\theta$ such that ${u\inps{z_0 + te^{i\theta}}} > {u\inps{z_0}}$ [resp., ${u\inps{z_0 + te^{i\theta}}} < {u\inps{z_0}}$] for sufficiently small positive $t$.  As before, let $C_a(u, z_0)$ and $C_d(u, z_0)$ represent the cones of, respectively, ascent and descent of $u$ at $z_0$. \par
In what follows, assume $\Omega$ to be an open, connected subset of $\C$.  We observe that $Re(f)$ and $e^f$ have the same cones of ascent and descent.
\begin{lemma} \label{lemma:uexp}
	Let $u:\Omega\rightarrow\R$ be the real part of a holomorphic function $f:\Omega\rightarrow\C$.  Then for $z_0\in\Omega$, $C_a(u, z_0) = C_a\inps{e^{f}, z_0}$ and $C_d(u, z_0) = C_d\inps{e^{f}, z_0}$.
	\begin{proof}
		Note that
		\begin{align}
			\aval{e^{f(z)}} &= \aval{e^{u(z) + iv(z)}} \\
			&= \aval{e^{u(z)}} \\
			&= e^{u(z)},
		\end{align}
		so that $\aval{e^{f(z_1)}} - \aval{e^{f(z_0)}}$ and $u(z_1) - u(z_0)$ have the same sign for all $z_0, z_1\in\C$.
	\end{proof}
\end{lemma}
The following two lemmas relate the behavior of $e^f$ and $f$ near $z_0$.
\begin{lemma} \label{lemma:expd}
	If $k\geq 1$ and $f(z)$ is $k$ times differentiable, then
	\begin{equation}
		\left[\exp\inps{f(z)}\right]^{(k)} = \sum_{j=0}^{k-1}\binom{k-1}{j} f^{\inps{k-j}}\inps{z}\cdot \left[\exp\inps{f(z)}\right]^{\inps{j}}. \label{eq:sum}
	\end{equation}
	\begin{proof}
	For $k = 1$, we simply apply the Chain Rule.  For $k > 1$, we may apply the General Leibniz Rule to $\left[\exp\inps{f(z)}\right]'$:
	\begin{align}
		\left[\exp\inps{f(z)}\right]^{(k)} &= \left[f'(z)\cdot \exp\inps{f(z)}\right]^{(k-1)} \\
		 &= \sum_{j=0}^{k-1}\binom{k-1}{j}\left[f'(z)\right]^{\inps{k-1-j}} \left[\exp\inps{f(z)}\right]^{(j)} \\
		 &= \sum_{j=0}^{k-1}\binom{k-1}{j} f^{\inps{k-j}}\inps{z}\cdot \left[\exp\inps{f(z)}\right]^{\inps{j}}. 
	\end{align}
	\end{proof}
\end{lemma}
\begin{lemma} \label{lemma:exp}
	Let $f:\Omega\rightarrow\C$ be holomorphic at $z_0\in\Omega$ and define $g(z):=\exp\inps{f(z)}$.  Suppose that, for some $k\geq 1$ and $z_0\in\Omega$,
\begin{enumerate}[label=(\arabic*)]
	\item $f^{(k)}(z_0)\neq 0$ and
	\item $f^{(\ell)}(z_0) = 0$ for each $\ell$ such that $1 \leq \ell < k$.
\end{enumerate}
Then
\begin{enumerate}[label=(\roman*)]
	\item $g^{(k)}(z_0) = g(z_0)\cdot f^{(k)}(z_0)\neq 0$ \label{item:kder} and
	\item $g^{(\ell)}(z_0) = 0$ for each $\ell$ such that $1 \leq \ell < k$. \label{item:zeroders}
\end{enumerate}
	\begin{proof}
		Both items follow readily from applying the hypotheses to Equation \eqref{eq:sum}.  
	\end{proof}
\end{lemma}

Using the preceding lemmas, we extend the Geometric Modulus Principle to harmonic functions:
\begin{theorem}[Geometric Modulus Principle for harmonic functions] \label{thm:hgmp}
	Let $u:\Omega\rightarrow\R$ be the real part of a holomorphic function $f:\Omega\rightarrow\C$ such that $k:=\min_{\ell \geq 1} \sbld{\ell : f^{(\ell)}(z_0) \neq 0}$.  Then $C_a(u, z_0)$ and $C_d(u, z_0)$ partition the unit circle into $2k$ alternating sectors of angle $\frac{\pi}{k}$.  Specifically, if $se^{i\beta}$ is the polar form of $f^{(k)}(z_0)$ then $\theta$ is an argument of ascent if
		\begin{equation} \label{eq:aor}
			\frac{-2\beta - \pi}{2k} < \theta < \frac{-2\beta + \pi}{2k}\quad \inps{\mbox{mod }\frac{2\pi}{k}}.
		\end{equation}
The interiors of the remaining sectors are arguments of descent.
\begin{proof}
	By Lemma \ref{lemma:uexp}, it suffices to apply Theorem \ref{thm:egmp} to $e^f$.  Combining Lemma \ref{lemma:exp} with Theorem \ref{thm:egmp}, we see that $e^f =: g$ yields $2k$ equal sectors whose angle of rotation is determined by the argument of
	\begin{align}
		\overline{g(z_0)}g^{(k)}(z_0) &= \overline{g(z_0)} g(z_0) f^{(k)}(z_0) \\
		&= \aval{g(z_0)}^2\cdot f^{(k)}(z_0).
	\end{align}
	Since we need only the argument of this quantity, we may drop the real coefficient $\aval{g(z_0)}^2$.
\end{proof}
\end{theorem}
\paragraph{Remark} Theorem \ref{thm:hgmp} appears (in a non-constructive form) as Theorem 12.1 of \cite{pvf}, where it is proved by properties of conformal maps. \par \vspace{1em}

\section{Corollaries} \label{sec:co}
In his original paper, Kalantari uses the Geometric Modulus Principle to give a two-line proof of the Maximum Modulus Principle for polynomials.  Having expanded his theorem to apply to all holomorphic functions, we generalize his proof\footnote{\cite[Theorem 2]{geomod}} to a proof of the general Maximum Modulus Principle.

\begin{theorem}[Maximum Modulus Principle]
	Let $\Omega$ be a connected, open subset of $\C$; let $f:\Omega\rightarrow\C$ be holomorphic.  If
		$\aval{f(z_0)} \geq \aval{f(z)}$
	for all $z$ in a neighborhood of some $z_0\in \Omega$, then $f$ is constant on $\Omega$.
	\begin{proof}
		By Theorem \ref{thm:egmp}, if $f$ is non-constant then it has at least one direction of ascent at each $z_0\in \Omega$; in particular, $\aval{f}$ cannot attain a local maximum in $\Omega$.
	\end{proof}
\end{theorem}

By analogous reasoning, we prove the Maximum Principle for harmonic functions.
\begin{theorem}[Maximum Principle]
	If the harmonic function $u:\Omega\rightarrow\R$ attains a local extremum at $z_0\in\Omega$, then $u$ is constant.
	\begin{proof}
		If $u$ is non-constant, then it is locally the real part of a non-constant holomorphic function.  Thus, by Theorem \ref{thm:hgmp}, it has at least one argument of ascent and at least one argument of descent.
	\end{proof}
\end{theorem}

Another straightforward corollary of Theorem \ref{thm:hgmp} is that the zeros of non-constant harmonic functions are never isolated, in contrast to the holomorphic case.
\begin{theorem}
	If the non-constant harmonic function $u:\Omega\rightarrow\C$ has a zero at $z = z_0$, then it also has a zero in every punctured disc $\sbld{z : 0 < \aval{z - z_0} < t}$.
	\begin{proof}
		By Theorem \ref{thm:hgmp}, the cones of ascent and descent of $u$ at $z_0$ partition the unit circle into at least two sectors.  By continuity, these sectors are bounded by the level curve $\sbld{z\in\Omega : u(z) = 0}$.
	\end{proof}
\end{theorem}

\section{Conclusion}
By expanding the Geometric Modulus Principle, we obtain stronger versions of two of the most powerful theorems in analysis: the Maximum Modulus Principle for holomorphic functions and the Maximum Principle for harmonic functions.  Given the implications of these theorems, of which Liouville's Theorem and the Schwarz lemma stand out in particular, we expect that careful application of the Geometric Modulus Principle will yield more results like those of Section \ref{sec:co}, or allow for briefer or more constructive proofs of existing results.  

\section*{Acknowledgements}
I would like to thank Yuefei Wang for his comments on an early draft.
\section*{Declarations}
\subsection*{Ethics approval and consent to participate}
Not applicable.
\subsection*{Consent for publication}
Not applicable.
\subsection*{Availability of data and materials}
Not applicable.
\subsection*{Competing interests}
None.
\subsection*{Funding}
For the 2020-21 academic year, I am receiving a fellowship from the SAS Department at Rutgers University.

\setcitestyle{numbers}
\bibliography{GeoMod}

\end{document}